\documentclass[10pt,leqno]{amsart}
\usepackage{amscd}
\usepackage{amssymb}
\usepackage{amsfonts}
\usepackage{latexsym}
\usepackage{verbatim}

\theoremstyle{plain}
\newtheorem{theorem}{Theorem}[section]
\newtheorem{definition}[theorem]{Definition}
\newtheorem{lemma}[theorem]{Lemma}
\newtheorem{prop}[theorem]{Proposition}
\newtheorem{cor}[theorem]{Corollary}
\newtheorem{rem}[theorem]{Remark}
\newtheorem{ex}[theorem]{Example}
\begin{document}
\title[Contact Calabi-Yau manifolds]{Contact Calabi-Yau manifolds and Special Legendrian submanifolds}
\author{Adriano Tomassini and Luigi Vezzoni}
\date{\today}
\address{Dipartimento di Matematica\\ Universit\`a di Parma\\ Viale G.
  P. Usberti 53/A\\
43100 Parma\\ Italy}
\email{adriano.tomassini@unipr.it}
\address{Dipartimento di Matematica \\ Universit\`a di Torino\\
Via Carlo Alberto 10\\
10123 Torino\\ Italy} \email{luigi.vezzoni@unito.it}
\thanks{This work was supported by the Projects M.I.U.R. ``Geometric Properties of Real and Complex Manifolds'', ``Riemannian
Metrics and Differentiable Manifolds'' and by G.N.S.A.G.A.
of I.N.d.A.M.}
\begin{abstract}
We consider a generalization of Calabi-Yau structures in the context
of Sasakian manifolds. We study deformations of a special
class of Legendrian submanifolds and classify invariant contact
Calabi-Yau structures on 5-dimensional nilmanifolds. Finally we
generalize to codimension $r$.
\end{abstract}
\maketitle
\newcommand\C{{\mathbb C}}
\newcommand\R{{\mathbb R}}
\newcommand\Z{{\mathbb Z}}
\newcommand\T{{\mathbb T}}
\newcommand{\de}[2]{\frac{\partial #1}{\partial #2}}
\newcommand\w{{\widetilde{J}}}
\newcommand{\ov}[1]{\overline{ #1}}
\newcommand{\Tk}{\mathcal{T}_{\kappa}}
\section{Introduction}
In their celebrated paper \cite{HL} Harvey and Lawson introduced the
concept of calibration and calibrated geometry. Namely, a
\emph{calibration} on an $n$-dimensional oriented Riemannian manifold $(M,g)$
is a closed $r$-form $\phi$ such that for any $x\in M$
$$
\phi_{x|V}\leq \mbox{Vol}(V)\,,
$$
where $V$ is an arbitrary oriented $r$-plane in $T_xM$. An oriented
submanifold $p\colon L\hookrightarrow M$ is said to be
\emph{calibrated} by $\phi$ if $p^*(\phi)=$Vol$(L)$. Compact
calibrated submanifolds have the important property of minimizing
volume in their homology class. As a typical example, the real part
of holomorphic volume form of a Calabi-Yau manifold is a
calibration; the corresponding calibrated submanifolds are said to
be \emph{special Lagrangian}. In \cite{McL} McLean studied special
Lagrangian submanifolds (and other special calibrated geometries)
showing that the Moduli space of deformations of special Lagrangian
manifolds of a fixed compact one $L$ is a smooth manifold of
dimension equal to the first Betti number of $L$. \\

In this paper we consider a generalization of Calabi-Yau structures
in the context of Sasakian manifolds. Recall that a \emph{Sasakian
structure} on a $2n+1$-dimensional manifold $M$ is a pair
$(\alpha,J)$, where $\alpha$ is a contact form on $M$ and  $J$ is an
integrable complex structure on
$\xi=\ker\alpha$ calibrated by $\kappa=\frac{1}{2}d\alpha$. This is equivalent to require the following data:
a quadruple $(\alpha,g,R,J)$, where $\alpha$ is a $1$-form, $g$ is a
Riemannian metric, $R$ is a unitary Killing vector field,
$J\in\mbox{End}(TM)$ satisfying
$$
J^2=-Id+\alpha\otimes R\,,\quad
g(J\,\cdot,J\,\cdot)=g(\cdot,\cdot)-\alpha\otimes\alpha\,,\quad
\alpha(R)=1
$$
and such that the metric cone $(M\times\R^+,r^2g+dr\otimes dr)$
endowed with the almost complex structure
$\widetilde{J}=J-r\alpha\otimes\partial_r+(1/r) dr\otimes R$ is
K\"ahler, where we extend $J$ by $J(\partial_r)=0$ (see e.g.
\cite{BL}, \cite{BG}, \cite{MSY}). These manifolds have been
studied by many authors (see e.g. \cite{BL}, \cite{BGM},
\cite{FOW}, \cite{MSY1}, \cite{MSY} and the references included).

We consider contact Calabi-Yau manifolds which are a special class of
Sasakian manifolds: namely a \emph{contact Calabi-Yau manifold} is a
$2n+1$-dimensional Sasakian manifold $(M,\alpha,J)$ endowed with a
closed basic complex volume form $\epsilon$. It turns out that these manifolds
are a special class of null-Sasakian $\alpha$-Einstein manifolds. As a direct consequence of the above definition,
in a contact Calabi-Yau manifold $(M,\alpha,J,\epsilon)$ the real
part of $\epsilon$ is a calibration. Furthermore, we have that
an $n$-dimensional submanifold $p\colon L\hookrightarrow M$ of a
contact Calabi-Yau manifold admits an orientation making it a calibrated submanifold
by $\Re\mathfrak{e}\,\epsilon$ if and only if
$$
p^*(\alpha)=0\,,\quad p^*(\Im\mathfrak{m}\,\epsilon)=0\,.
$$
In such a case $L$ is said to be a \emph{special Legendrian
  submanifold}.
We prove that:
\vskip0.3cm
\noindent{\it The Moduli space of deformations of special
  Legendrian submanifolds near a fixed compact one $L$ is a smooth
\vspace{0.3cm}
  $1$-dimensional manifold}. \\
Moreover we get the following extension theorem: \vskip0.3cm
\noindent{\it Let $(M,\alpha_t,J_t,\epsilon_t)$ be a smooth family of
contact Calabi-Yau manifolds and let $p\colon L\hookrightarrow
(M,\alpha_0,J_0,\epsilon_0)$ be a compact special Legendrian
submanifold. Then there exists a smooth family of special Legendrian
submanifolds $p_t\colon L\hookrightarrow
(M,\alpha_t,J_t,\epsilon_t)$ that extends $p\colon L\hookrightarrow M$ if and only if the
cohomology class \vspace{0.3cm}
$[p^*(\Im\mathfrak{m}\,\epsilon)]$ vanishes.}\\
This can be considered a contact version of a theorem of Lu Peng
(see \cite{LP})
in
Calabi-Yau manifolds (see also \cite{Pa}). \\

In section 2 we fix some notation on contact and Sasakian geometry.
In section 3 we define contact Calabi-Yau manifolds and we obtain
some simple topological obstructions to the existence of contact
Calabi-Yau structures on odd-dimensional manifolds. As a corollary,
we get that there are no contact Calabi-Yau structures on
odd-dimensional spheres. In section 4 we study the Moduli space of
special Legendrian submanifolds, proving the theorems stated
above.\\ In section 5 we classify the 5-dimensional nilmanifolds
carrying an invariant contact Calabi-Yau
structure. The proof is based on theorems 21 and 23 of \cite{CS}. \\
In the last section we generalize the previous definition to the
case of codimension $r$ proving an extension theorem.\\
Some
examples of contact Calabi-Yau manifolds and special Legendrian
submanifolds are carefully described.
\section{Preliminaries}
Let $M$ be a manifold of dimension $2n+1$. A \emph{contact
structure} on $M$ is a distribution $\xi\subset
TM$ of dimension $2n$, such that the defining 1-form $\alpha$ satisfies
\begin{equation}
\label{contact}
\alpha\wedge(d\alpha)^n\neq 0\,.
\end{equation}
A 1-form $\alpha$ satisfying \eqref{contact} is said to be a
\emph{contact form} on $M$. Let $\alpha$ be a contact form on $M$;
then there exists a unique vector field $R_{\alpha}$ on $M$ such that
$$
\alpha(R_{\alpha})=1\,,\quad \iota_{R_{\alpha}}d\alpha=0\,,
$$
where $\iota_{R_{\alpha}}d\alpha$ denotes the contraction of $d\alpha$
along $R_{\alpha}$. By definition $R_{\alpha}$ is called the \emph{Reeb vector
field} of the contact form $\alpha$.\\
A \emph{contact manifold} is a pair $(M,\xi)$ where $M$ is a
$2n+1$-dimensional manifold and $\xi$ is a contact structure. Let
$(M,\xi)$ be a contact manifold and fix a defining (contact) form
$\alpha$. Then the 2-form $\kappa=\frac{1}{2}d\alpha$ defines a symplectic form
on the contact structure $\xi$; therefore the pair $(\xi,\kappa)$
is a symplectic vector bundle over $M$. A \emph{complex structure}
on $\xi$ is the datum of $J\in\mbox{End}(\xi)$ such that $J^2=-I_{\xi}$.
\begin{definition}
Let $\alpha$ be a contact form on $M$, with $\xi=\ker\alpha$ and let
$\kappa=\frac{1}{2}d\alpha$. A complex structure $J$ on $\xi$ is said to be
$\kappa$-\emph{calibrated} if
$$
g_J[x](\cdot,\cdot):=\kappa[x](\cdot,J_x\cdot)
$$
is a $J_x$-Hermitian inner product on $\xi_{x}$ for any $x\in M$.
\end{definition}
The set of
$\kappa$-calibrated complex structures on $\xi$ will be denoted by $\mathfrak{C}_{\alpha}(M)$.\\
If $J$ is a complex structure on $\xi=\ker\alpha$, then we extend it
to an endomorphism of $TM$ by setting
$$
J(R_{\alpha})=0\,.
$$
Note that such a $J$  satisfies
$$
J^2=-I+\alpha\otimes R_{\alpha}\,.
$$
If $J$ is $\kappa$-calibrated, then it induces a Riemannian
metric $g$ on $M$ given by
\begin{equation}\label{g}
g:=g_J+\alpha\otimes\alpha\,.
\end{equation}
Furthermore the Nijenhuis tensor of $J$ is defined by
$$
N_{J}(X,Y)=[JX,JY]-J[X,JY]-J[Y,JX]+J^2[X,Y]
$$
for any $X,Y\in TM$. We recall the following
\begin{definition}
A Sasakian structure on a $2n+1$-dimensional manifold $M$ is a pair
$(\alpha,J)$, where
\begin{itemize}
\item $\alpha$ is a contact form;
\vskip0.1cm
\item $J\in\mathfrak{C}_{\alpha}(M)$ satisfies $N_{J}=-d\alpha\otimes R_{\alpha}$.
\end{itemize}
The triple $(M,\alpha,J)$ is said to be a \emph{Sasakian manifold}.
\end{definition}
\noindent For other characterizations of Sasakian structure see e.g.
\cite{BL} and \cite{BG}.

We recall now the definition of basic $r$-forms.
\begin{definition}
Let $(M,\xi)$ be a contact manifold. A differential $r$-form
$\gamma$ on $M$ is said to be \emph{basic} if
$$
\iota_{R_{\alpha}}\gamma=0\,,\quad \mathcal{L}_{R_{\alpha}}\gamma=0\,,
$$
where $\mathcal{L}$ denotes the Lie derivative and $R_{\alpha}$ is
the Reeb vector field of an arbitrary contact form defining $\xi$.
\end{definition}
\noindent We will denote by $\Lambda^{r}_{B}(M)$ the set of basic
$r$-forms on $(M,\xi)$. Note that
$$
d\,\Lambda^r_B(M)\subset\Lambda^{r+1}_{B}(M)\,.
$$
The cohomology $H_{B}^{\bullet}(M)$ of this complex is called the \emph{basic cohomology}
of $(M,\xi)$.

If $(M,\alpha,J)$ is a Sasakian manifold, then
$$
J(\Lambda^r_B(M))=\Lambda^r_{B}(M)\,,
$$
where, as usual, the action of $J$ on $r$-forms is defined by
$$
J\phi(X_{1},\dots,X_{r})=\phi(JX_1,\dots,JX_r)\,.
$$
Consequently $\Lambda^r_B(M)\otimes\C$ splits as
$$
\Lambda^r_B(M)\otimes\C=\bigoplus_{p+q=r}\Lambda^{p,q}_J(\xi)
$$
and, according with this gradation, it is possible to define the
cohomology groups $H^{p,q}_{B}(M)$. The $r$-forms belonging to
$\Lambda^{p,q}_J(\xi)$ are said to be of \emph{type} $(p,q)$
\emph{with respect to} $J$. Note that
$\kappa=\frac{1}{2}d\alpha\in\Lambda^{1,1}_J(\xi)$ and it determines a
non-vanishing cohomology class in $H_{B}^{1,1}(M)$.
The Sasakian structure $(\alpha,J)$ also induces a natural connection $\nabla^{\xi}$ on
$\xi$ given by
$$
\nabla^{\xi}_XY=\begin{cases} (\nabla_{X}Y)^{\xi}\mbox{ if }X\in\xi\\
[R_{\alpha},Y]\,\,\,\mbox{ if }X=R_{\alpha}\,,
\end{cases}
$$
where the subscript $\xi$ denotes the projection onto $\xi$. One
easily gets
$$
\nabla^{\xi}_{X}J=0\,,\quad
\nabla^{\xi}_Xg_J=0\,,\quad\nabla^{\xi}_Xd\alpha=0\,,\quad \nabla^{\xi}_XY-\nabla^{\xi}_{Y}X=[X,Y]^{\xi}\,,
$$
for any $X,Y\in TM$. Consequently we have
$$
\mbox{Hol}(\nabla^{\xi})\subseteq \mbox{U}(n).
$$
Moreover the \emph{transverse Ricci
tensor} $\mbox{Ric}^T$ is defined as
$$
\mbox{Ric}^T(X,Y)=\sum_{i=1}^{2n}g(\nabla^{\xi}_X\nabla^{\xi}_{e_i}e_i-\nabla^{\xi}_{e_i}\nabla^{\xi}_{X}e_i-\nabla^{\xi}_{[X,e_i]}e_i,Y)\,,
$$
for any $X,Y\in\xi$, where $\{e_1,\dots,e_{2n}\}$ is an arbitrary orthonormal frame of $\xi$. It
is known that $\mbox{Ric}^T$ satisfies
$$
\mbox{Ric}^T(X,Y)=\mbox{Ric}(X,Y)+2\,g(X,Y)\,,
$$
for any $X,Y\in\xi$, where Ric denotes the Ricci tensor
of the Riemannian metric $g=g_J+\alpha\otimes\alpha$. Let us denote
by $\rho^T$ the Ricci form of $\mbox{Ric}^T$, i.e.
$$
\rho^T(X,Y)=\mbox{Ric}^T(JX,Y)=\mbox{Ric}(JX,Y)+2\,\kappa(X,Y)\,,
$$
for any $X,Y\in\xi$. We recall that $\rho^T$ is a closed form such that $\frac{1}{2\pi}\rho$ represents
the first Chern class of $(\xi,J)$
(see e.g.
\cite{EKA}); this form is called the \emph{transverse Ricci form} of $(\alpha,J)$.
\begin{definition}
The basic cohomology class
$$
c_{1}^B(M)=\frac{1}{2\pi}[\rho^T]\in H^{1,1}_B(M)
$$
is called the \emph{first basic Chern class} of $(M,\alpha,J)$ and,
if it vanishes, then $(M,\alpha,J)$ is said to be \emph{null-Sasakian}.
\end{definition}
Furthermore we recall that a Sasakian manifold is called
\emph{$\alpha$-Einstein} if there exist $\lambda,\nu\in
C^{\infty}(M,\R)$ such that
$$
\mbox{Ric}=\lambda g+\nu \alpha\otimes \alpha\,.
$$
For general references on these topics see e.g. \cite{BGN} and \cite{BGM}.

Finally, recall
that a submanifold $p\colon L\hookrightarrow M$
of a $2n+1$-dimensional contact manifold
$(M,\xi)$ is said to be \emph{Legendrian} if :
\begin{enumerate}
\item[1)] dim$_{\R}\,L=n$,
\vspace{0.2 cm}
\item[2)]$p_*(TL)\subset\xi$.
\end{enumerate}
\vskip0.1truecm
Observe that, if $\alpha$ is a defining form of the contact structure
$\xi$, then condition 2) is equivalent to say that $p^*(\alpha)=0$.\\
Hence Legendrian  submanifolds are the analogue of Lagrangian
submanifolds in contact geometry.
\section{Contact Calabi-Yau manifolds}
In this section we study contact Calabi-Yau manifolds. As
already explained in the introduction, these manifolds are a natural
generalization of the Calabi-Yau ones in the context of contact
geometry. Roughly speaking a  contact Calabi-Yau manifold is a Sasakian
manifold endowed with a basic closed complex volume
form. We can give now the following
\begin{definition}
A \emph{contact Calabi-Yau manifold} is a quadruple
$(M,\alpha,J,\epsilon)$, where
\begin{itemize}
\item $(M,\alpha,J)$ is a $2n+1$-dimensional Sasakian  manifold;
\vspace{0.1cm}
\item $\epsilon\in\Lambda^{n,0}_{J}(\xi)$ is a nowhere vanishing basic form on
  $\xi=\ker\alpha$ such that
$$
\begin{cases}
\epsilon\wedge\overline{\epsilon}=c_{n}\,\kappa^n\\
d\epsilon=0\,,
\end{cases}
$$
where $c_{n}=(-1)^{\frac{n(n+1)}{2}}(2i)^n$ and $\kappa=\frac{1}{2}d\alpha$.
\end{itemize}
\end{definition}
Now we will describe a couple of examples.
\begin{ex}\emph{
Consider $\R^{2n+1}$ endowed with the standard Euclidean coordinates
$\{x_1,\dots,x_n,y_1,\dots,y_n,t\}$. Let
$$
\alpha_{0}=2dt-2\sum_{i=1}^n y_i\,dx_i\,.
$$
be the \emph{standard contact form} on $\R^{2n+1}$ and let
$\xi_{0}=\ker\alpha_0$. Then $\xi_{0}$ is spanned by
$$
\{y_1\,\partial_t+\partial_{x_{1}},\,
\dots,\,y_n\,\partial_t+\partial_{x_{n}},\,\partial_{y_{1}},\,\dots,\,\partial_{y_{n}}\}\,.$$
For simplicity, set $V_i=y_i\,\partial_t+\partial_{x_{i}}$,
$W_j=\partial_{y_{j}}$, $i,j=1,\dots,n$ and
$$
\begin{cases}
J_0(V_{r})=W_{r}\\
J_0(W_{r})=-V_r
\end{cases}
\quad r=1,\dots n\,.
$$
Then $J_0$ defines a complex structure in
$\mathfrak{C}_{\alpha}(M)$. Since the space of transverse 1-forms is
spanned by $\{dx_1,\,\dots,\,dx_n,\,dy_1,\dots,\,dy_n\}$, then the
complex valued form
$$
\epsilon_0:=(dx_1+idy_1)\wedge\dots\wedge(dx_n+idy_n)
$$
is of type $(n,0)$ with respect to $J_0$ and it satisfies
$$
\begin{cases}
\epsilon_0\wedge\overline{\epsilon}_0=c_{n}\,\kappa_{0}^{n}\\
d\epsilon_0=0\,,
\end{cases}
$$
where $\kappa_0=\frac{1}{2}d\alpha_0$.
Therefore $(\R^{2n+1},\alpha_{0},J_{0},\epsilon_{0})$ is a contact
Calabi-Yau manifold.}
\end{ex}
The following will describe a compact contact Calabi-Yau manifold.
\begin{ex}\label{exam}
\emph{
Let
$$
H(3):=\left\{
A=\left(
\begin{array}{ccc}
1  &x &y\\
0  &1 &z\\
0  &0 &1 \\
\end{array}
\right)\,\,|\,\,x,y,z\in \R
\right\}\
$$
be the 3-dimensional Heisenberg group and let $M=H(3)/\Gamma$, where $\Gamma$
denotes the subgroup of $H(3)$ given by the matrices
with integral entries. The 1-forms
$\alpha_1=dx$,
$\alpha_2=dy$,
$\alpha_3=x\,dy-dz$ are $H(3)$-invariant and therefore they define a
global coframe on $M$. Then $\alpha=2\alpha_3$ is a contact form whose contact
distribution $\xi$ is spanned by $V=\partial_{x}$,
$W=\partial_y+x\,\partial_z$. Again
$$
\begin{cases}
J(V)=W\\
J(W)=-V
\end{cases}
$$
defines a $\kappa$-calibrated complex structure on $\xi$ and
$\epsilon=\alpha_1+i\,\alpha_2$ is a $(1,0)$-form on $\xi$ such that
$(M,\alpha,J,\epsilon)$ is a contact Calabi-Yau manifold.}
\end{ex}
The last example gives an invariant contact Calabi-Yau structure on
a nilmanifold. It can be generalized to the dimension $2n+1$ in this
way: let $\mathfrak{g}$ be the Lie algebra spanned by
$\{X_{1},\dots,X_{2n+1}\}$ with
$$
[X_{2k-1},X_{2k}]=-X_{2n+1}
$$
for $k=1,\dots, n$ and the other brackets are zero. Then $\mathfrak{g}$
is a $2n+1$-dimensional nilpotent Lie algebra with rational constant
structures and, by Malcev theorem, it follows that
if $G$ is the simply connected Lie
group with Lie algebra $\mathfrak{g}$, then $G$ has a compact
quotient. Let $\{\alpha_1,\dots,\alpha_{2n+1}\}$ be the dual basis of
$\{X_1,\dots,X_{2n+1}\}$.
Then we immediately get
$$
d\alpha_1=0\,,\dots,\, d\alpha_{2n}=0\,,\quad d\alpha_{2n+1}=\sum_{k=1}^{n}\alpha_{2k-1}\wedge\alpha_{2k}\,.
$$
Hence
$$
\alpha=2\,\alpha_{2n+1}\,,
$$
the endomorphism $J$ of $\xi=\ker\alpha$ defined by
$$
\begin{cases}
J(X_{2k-1})=X_{2k}\\
J(X_{2k})=-X_{2k-1}
\end{cases}
$$
for $k=1,\dots,n$ and the complex form
$$
\epsilon=(\alpha_1+i\,\alpha_2)\wedge\dots\wedge (\alpha_{2n-1}+i\,\alpha_{2n})
$$
define a contact Calabi-Yau structure on any compact nilmanifold
associated with $\mathfrak{g}$.\\

The following proposition gives simple topological obstructions in
order that a compact $2n+1$-dimensional manifold $M$ carries a contact
Calabi-Yau structure.
\begin{prop}\label{ostruzioni}
Let $M$ be a $2n+1$-dimensional compact manifold. Assume that $M$
admits a contact Calabi-Yau structure; then the following hold
\begin{enumerate}
\item[1.] if $n$ is \emph{even}, then $b_{n+1}(M)> 0$;
\vskip0.1cm
\item[2.] if $n$ is \emph{odd},  then
$$
\begin{cases}
b_{n}(M)\geq 2\\
b_{n+1}(M)\geq 2\,,
\end{cases}
$$
\end{enumerate}
where $b_{j}(M)$ denotes the $j^{th}$ Betti number of $M$.
\end{prop}
\begin{proof}
Let $(\alpha,J,\epsilon)$ be a contact Calabi-Yau structure on
$M$ and let $\xi=\ker\alpha$. Set $\Omega=\Re\mathfrak{e}\,\epsilon$;
then, since $\epsilon\in\Lambda^{n,0}_J(\xi)$ , we have
$\epsilon=\Omega+i\,J\Omega$.
In view of the assumption $d\epsilon=0$, we obtain $d\Omega=dJ\Omega=0$ and since
$d\alpha\in\Lambda^{1,1}_{J}(\xi)$ it follows that
$$
\Omega\wedge d\alpha=J\Omega\wedge d\alpha=0\,.
$$
Hence
$$
d(\Omega\wedge\alpha)=d(J\Omega\wedge \alpha)=0\,.
$$
Furthermore we have
$$
\begin{array}{lll}
&\epsilon\wedge\overline{\epsilon}=\Omega\wedge\Omega+J\Omega\wedge
J\Omega \quad &\mbox{ if $n$ is even;}\\[5pt]
&\epsilon\wedge\overline{\epsilon}=-2i\,\Omega\wedge J\Omega &\mbox{ if
  $n$ is odd}.
\end{array}
$$
\begin{enumerate}
\item[1.] If  $n$ is even, then
  $\alpha\wedge(\Omega\wedge\Omega+J\Omega\wedge J\Omega)$ is a volume form on $M$. Assume that
  the cohomology classes $[\Omega\wedge\alpha]$, $[J\Omega\wedge\alpha]$ vanish; then
  there exist $\beta,\gamma\in \Lambda^{n}(M)$ such that
$$
\alpha\wedge\Omega=d\beta\,,\quad \alpha\wedge J\Omega=d\gamma\,.
$$
By Stokes theorem we have
$$
\begin{aligned}
0\neq\int_{M}\alpha\wedge\Omega\wedge\Omega+\alpha\wedge J\Omega\wedge
J\Omega
&=\int_{M}d\beta\wedge\Omega+d\gamma\wedge J\Omega\\
&=\int_{M}d(\beta\wedge\Omega)+d(\gamma\wedge J\Omega)=0\,,
\end{aligned}
$$
which is absurd. Therefore one of $[\Omega\wedge\alpha]$,
$[J\Omega\wedge\alpha]$ does not vanish.
Consequently $b_{n+1}(M)> 0$.
\vspace{0.3cm}
\item[2.] Let $n$ be odd. We prove that the cohomology classes
  $[\Omega]$ and $[J\Omega]$ are $\R$-independent. Assume that there
  exist $a,b\in\R$ such that $a[\Omega]+b[J\Omega]=0$, $(a,b)\neq (0,0)$. Then
  there exists
  $\beta\in\Lambda^{n-1}(M)$ such that
$$
a\,\Omega+b\,J\Omega=d\beta\,.
$$
We may assume that $a=1$, so that $\Omega=d\beta-bJ\Omega$. Stokes
theorem implies
$$
0\neq \int_{M}\alpha\wedge\Omega\wedge J\Omega
=\int_{M}\alpha\wedge d\beta\wedge J\Omega=
-\int_{M}d(\alpha\wedge\beta\wedge J\Omega)=0
$$
which is a contradiction. Hence $b_n(M)\geq 2$.  With the same argument, it is possible to prove
that $b_{n+1}(M)\geq 2$ by showing that $[\Omega\wedge\alpha]$ and
$[J\Omega\wedge\alpha]$ are $\R$-independent in $H^{n+1}(M,\R)$.
\end{enumerate}
\end{proof}
The following is an immediate consequence of proposition \ref{ostruzioni}.
\begin{cor}
A $3$-dimensional compact manifold $M$ admitting contact Calabi-Yau
structure has $b_{1}(M)\geq 2$. In particular, there are no
compact $3$-dimensional simply connected contact Calabi-Yau manifolds.\\
Moreover, the $2n+1$-dimensional sphere has no contact Calabi-Yau structures.
\end{cor}
The following proposition implies that the transverse Ricci tensor of a contact Calabi-Yau manifold vanishes
\begin{prop}\label{copiajoice}
Let $(M,\alpha,J)$ be a $2n+1$-dimensional Sasakian manifold and $\xi=\ker\alpha$. The
following facts are equivalent:
\begin{enumerate}
\item[$1.$] $\mbox{\emph{Hol}}^0(\nabla^{\xi})\subseteq \mbox{\emph{SU}}(n)$
\vskip0.1cm
\item[$2.$] $\mbox{\emph{Ric}}^T=0$.
\end{enumerate}
\end{prop}
\begin{proof}
The connection $\nabla^{\xi}$ induces a connection $\nabla^{K}$ on
$\Lambda^{n,0}_J(\xi)$ which has Hol$(\nabla^K)\subseteq \mbox{U}(1)$. Since
Hol$^{0}(\nabla^K)$ and Hol$^0(\nabla^{\xi})$ are related by
$$
\mbox{Hol}^{0}(\nabla^K)=\det (\mbox{Hol}^0(\nabla^{\xi}))\,,
$$
where $\det$ is the map induced by the determinant
U$(n)\to\,$U$(1)$, then it follows that $\mbox{Hol}^0(\nabla^{\xi})\subseteq \mbox{SU}(n)$ if and only if
 $\mbox{Hol}^{0}(\nabla^K)=\{1\}$ and in this case
$\nabla^K$ is flat. As in the K\"ahler case it can be showed using transverse holomorphic coordinates (see e.g. \cite{EKA}, \cite{FOW})
 that the
curvature form of $\nabla^K$ coincides with the transverse Ricci form of
$(\alpha,J)$. Hence $\mbox{Hol}^0(\nabla^{\xi})\subseteq
\mbox{SU}(n)$ if and only if $\mbox{Ric}^T=0$.
\end{proof}
As a consequence of the last proposition we have the following
\begin{cor}
Let $(M,\alpha,J,\epsilon)$ be a contact Calabi-Yau manifold. Then $(M,\alpha,J)$ is null-Sasakian and the metric $g$ induced by $(\alpha,J)$
is $\alpha$-Einstein with $\lambda=-2$ and $\nu=2n+2$. In particular the scalar curvature of the metric $g$ associated to $(\alpha,J)$ is
equal to $-2n$.
\end{cor}
\section{Deformations of  special Legendrian submanifolds}
In this section we are going to study the geometry of Legendrian
submanifolds in a contact Calabi-Yau ambient. We will prove a
contact version of McLean and Lu Peng theorems (see \cite{McL} and
\cite{LP}).

Let $(M,\alpha,J,\epsilon)$ be a contact Calabi-Yau  manifold of
dimension $2n+1$.\\
It easy to see that for any oriented $n$-plane $V\subset T_{x}M$
$$
\Re\mathfrak{e}\,\epsilon_{|V}\leq \mbox{Vol}(V)\,,
$$
where $\mbox{Vol}(V)$ is computed with respect to the metric $g$
induced by $(\alpha,J)$ on $M$. Hence $\Re\mathfrak{e}\,\epsilon$ is
a calibration on $(M,g)$ (see \cite{HL}). We have the following
\begin{prop}
Let $p\colon L\hookrightarrow M$ be an $n$-dimensional submanifold.
The following facts are equivalent
\begin{itemize}
\item[1.] the submanifold satisfies
$$
\begin{cases}
p^*(\alpha)=0\\
p^*(\Im \mathfrak{m}\,\epsilon)=0\,,
\end{cases}
$$
\item[2.] there exists an orientation on $L$ making it calibrated by $\Re\mathfrak{e}\,\epsilon$.
\end{itemize}
\end{prop}
We can give the following
\begin{definition}
An $n$-dimensional submanifold $p\colon L\hookrightarrow M$ is said to
be \emph{special Legendrian} if
$$
\begin{cases}
p^*(\alpha)=0\\
p^*(\Im \mathfrak{m}\,\epsilon)=0\,.
\end{cases}
$$
\end{definition}
It follows that compact special Legendrian submanifolds minimize
volume in their homology class and that there are no compact special
Legendrian submanifolds in $(\R^{2n+1},\alpha_0,J_0,\epsilon_0)$.
\begin{ex}\label{SL}
\emph{
Let $(M=H(3)/\Gamma,\alpha,J,\epsilon)$ be the contact Calabi-Yau
manifold considered in the example \ref{exam}. Then
the submanifold
$$
L:=\left\{[A]\in M\,|\,
A=\left(
\begin{array}{ccc}
1  &x &0\\
0  &1 &0\\
0  &0 &1 \\
\end{array}
\right)
\right\}\simeq S^1
$$
is a compact special Legendrian submanifold.}
\end{ex}
Now we define the Moduli space of special Legendrian submanifolds.
\begin{definition}
Let $(M,\alpha,J,\epsilon)$ be a contact Calabi-Yau manifold and let
$p_0\colon L\hookrightarrow M$,
$p_1\colon L\hookrightarrow M$ be two special Legendrian submanifolds. Then
$p_1\colon L\hookrightarrow M$ is said to be a \emph{deformation} of $p_0\colon L\hookrightarrow M$ if there exists
a smooth map $F\colon L\times[0,1]\to M$ such that
\begin{itemize}
\item $F(\cdot,t)\colon L\times \{t\}\to M$ is a special Legendrian embedding
  for any $t\in[0,1]$;
\vspace{0.1cm}
\item $F(\cdot,0)=p_0$, $F(\cdot,1)=p_1$.
\end{itemize}
\end{definition}
Let $(M,\alpha,J,\epsilon)$ be a contact
Calabi-Yau manifold and let $p\colon L\hookrightarrow M$ be a fixed
compact special Legendrian submanifold. Set
$$
\begin{array}{l}
\mathfrak{M}(L):=\{\mbox{\rm special
Legendrian submanifolds of}\,\, (M,\alpha,J,\epsilon)\\[5pt]
\hskip1.6truecm\mbox{\rm which are deformations of}\,\, p\colon L \hookrightarrow M\}/\sim\,,
\end{array}
$$
where two embeddings are considered  equivalent if they differ by a
diffeomorphism of $L$; then by definition $\mathfrak{M}(L)$ is the
\emph{Moduli space of special Legendrian  submanifolds} which are
deformations of $p\colon L\hookrightarrow M$. We have the following
\begin{theorem}\label{McL}
Let $(M,\alpha,J,\epsilon)$ be a contact Calabi-Yau manifold and let
$p\colon L\hookrightarrow M$ be a compact special  Legendrian
submanifold. Then the Moduli space $\mathfrak{M}(L)$ is a
$1$-dimensional manifold.
\end{theorem}
The next lemma will be useful in the proof of theorem \ref{McL}:
\begin{lemma}[\cite{McL}, \cite{dB}]\label{linear}
Let $(V,\kappa)$ be a symplectic vector space and let $i\colon W\hookrightarrow V$
be a Lagrangian subspace. Then
\begin{enumerate}
\item[1.] $\tau\colon V/W\to W^*$ defined as
  $\tau([v])=i^*(\iota_v\kappa)$ is an isomorphism;
\item[2.] let $J$ be a $\kappa$-calibrated complex structure on $V$ and
  let $\epsilon\in\Lambda^{n,0}_{J}(V^*)$ satisfy
$$
i^*(\Im \mathfrak{m}\,\epsilon)=0\,,\quad \epsilon\wedge\ov{\epsilon}=c_n\frac{\kappa^n}{n!}\,.
$$
Then $\theta\colon V/W\to\Lambda^{n-1}(W^*)$ defined as
$\theta([v]):=i^*(\iota_v\Im\mathfrak{m}\,\epsilon)$ is an isomorphism. Moreover
for any $v\in V$, we have
$$
\theta([v])=-*\tau([v])\,,
$$
\end{enumerate}
where $*$ is computed with respect to
$i^*(g_{J}(\cdot,\cdot)):=i^*(\kappa(\cdot,J\cdot))$ and the volume form \emph{Vol}$(W):=i^*(\Re\mathfrak{e}\,\epsilon)$.
\end{lemma}
For the proof of lemma \ref{linear} we refer to \cite{McL} and \cite{dB}.
\begin{proof}[Proof of theorem $\ref{McL}$]
Let $\mathcal{N}(L)$ be the normal bundle to $L$. Then
$$
\mathcal{N}(L)=<R_{\alpha}>\oplus J (p_*(TL))\,,
$$
where $R_{\alpha}$ is the Reeb vector field of $\alpha$. Let $Z$ be
a vector field normal to $L$ and let $\exp_{Z}\colon L\to M$ be
defined as
$$
\exp_{Z}(x):=\exp_{x}(Z(x))\,.
$$
Let $U$ be a neighborhood of $0$ in $C^{2,\alpha}(<R_{\alpha}>)\oplus
C^{1,\alpha}(J (p_*(TL)))$
and let
$$
F\colon U\to C^{1,\alpha}(\Lambda^{1}(L))\oplus
C^{0,\alpha}(\Lambda^{n}(L))\,,
$$
be defined as
$$
F(Z)=(\exp_{Z}^*(\alpha),2\exp_{Z}^*(\Im \mathfrak{m}\,\epsilon))\,.
$$
We obviously have
$$
Z\in F^{-1}((0,0))\cap
C^{\infty}(\mathcal{N}(L))\iff\,\exp_{Z}(L)\mbox{
  is a special Legendrian submanifold.}
$$
Note that since $\exp_{Z}$ and $p$ are homotopic via $\exp_{t\,Z}$, we
have
$$
[\exp_{Z}^*(\Im\mathfrak{m}\,\epsilon)]=[p^*(\Im\mathfrak{m}\,\epsilon)]=0\,.
$$
Therefore
$$
F\colon U\to C^{1,\alpha}(\Lambda^{1}(L))\oplus dC^{1,\alpha}(\Lambda^{n-1}(L))\,.
$$
Let us compute the differential of the map $F$.
$$
F_{*}[0](Z)=\frac{d}{dt}(\exp_{t\,Z}^*(\alpha),2\exp_{t\,Z}^*(\Im\mathfrak{m}\,\epsilon))_{|t=0}=(p^*(\mathcal{L}_{Z}\alpha)
,2p^*(\mathcal{L}_{Z}\Im\mathfrak{m}\,\epsilon))\,,
$$
where $\mathcal{L}$ denotes the Lie derivative. We may write
$Z=JX+f\,R_{\alpha}$; then applying Cartan formula we obtain
$$
\begin{aligned}
F_{*}[0](Z)&=(p^*(\mathcal{L}_{Z}\alpha),2p^*(\mathcal{L}_Z\Im
\mathfrak{m}\,\epsilon))\\
&=(p^*(d\iota_Z\alpha+\iota_{Z}d\alpha),2p^*(d\iota_Z\Im
\mathfrak{m}\epsilon))\\
&=(p^*(d\iota_{JX+f\,R_{\alpha}}\alpha+\iota_{JX+f\,R_{\alpha}}d\alpha),2p^*(d\iota_{JX+f\,R_{\alpha}}\Im
\mathfrak{m}\,\epsilon))\\
&=(p^*(d\iota_{fR_{\alpha}}\alpha+\iota_{JX}d\alpha),2p^*(d\iota_{JX}\Im
\mathfrak{m}\,\epsilon))\\
&=(p^*(df+\iota_{JX}d\alpha),2dp^*(\iota_{JX}\Im
\mathfrak{m}\,\epsilon))\,.
\end{aligned}
$$
By applying lemma \ref{linear} we get
\begin{equation}
\label{nucleo}
F_{*}[0](Z)=(d(f\circ p)+p^*(\iota_{JX}d\alpha),-d*p^*(\iota_{JX}d\alpha))\,,
\end{equation}
where $*$ is the Hodge star operator with respect to the metric
$p^*(g_{J})$
and the volume form $p^*(\Re\mathfrak{e}\,\epsilon)$.\\
Now we show that $F_{*}[0]$ is surjective. Let $(\eta,d\gamma)\in
C^{1,\alpha}(\Lambda^{1}(L))\oplus dC^{1,\alpha}(\Lambda^{n-1}(L))$.
By the Hodge decomposition theorem we may assume
$$
d\gamma=-d*du\mbox{ with } u\in C^{3,\alpha}(L)
$$
and we have
$$
\eta=dv+d^*\beta+h(\eta)
$$
where $v\in C^{2,\alpha}(L)$, $\beta\in
C^{2,\alpha}(\Lambda^{2}(L))$ and $h(\eta)$ denotes the harmonic
component of $\eta$. Then we get
$$
\begin{aligned}
(\eta,d\gamma)=&(du-du+dv+d^*\beta+h(\eta),-d*du)\\
              =&(dv-du+du+d^*\beta+h(\eta),-d*(du+d^*\beta+h(\eta))\,.
\end{aligned}
$$
We can find  $f\in C^{2,\alpha}(p(L))$ and $X\in
C^{1,\alpha}(p_*(TL))$ such that
$$
\begin{aligned}
&f\circ p=v-u\\
&p^*(\iota_{JX}d\alpha)=du+d^*\beta+h(\eta)\,.
\end{aligned}
$$
Hence
$$
(\eta,d\gamma)=(d(f\circ p)+p^*(\iota_{JX}d\alpha),-d*p^*(\iota_{JX}d\alpha))
$$
and $F_{*}[0]$ is surjective. Therefore $(0,0)$ is a regular value of $F$.
Now we compute $\ker F_{*}[0]$. Formula \eqref{nucleo} implies that $Z\in\ker F_{*}[0]$ if and
only if
\begin{eqnarray}
&&\label{arr1} d(f\circ p)+p^*(\iota_{JX}d\alpha)=0\\
&&\label{arr2} d^*p^*(\iota_{JX}d\alpha)=0\,.
\end{eqnarray}
By applying $d^*$ to both sides of  \eqref{arr1} and taking
into account \eqref{arr2}  we get
$$
0=d^*d(f\circ p)+d^*p^*(\iota_{JX}d\alpha)=d^*d(f\circ p)\,,
$$
i.e.
$$
\Delta(f\circ p)=0\,.
$$
Since $L$ is compact $f$ is  constant. Hence
\eqref{arr1} reduces to
\begin{equation}\label{JX=0}
p^*(\iota_{JX}d\alpha)=0\,.
\end{equation}
The map
$$
\Theta\colon p_*(TL)\to \Lambda^1 (L)
$$
defined by
$$
\Theta(X)=p^*(\iota_{JX}d\alpha)
$$
is an isomorphism; hence equation \eqref{JX=0} implies $X=0$. Therefore
$Z=W+f\,R_{\alpha}$ belongs to $\ker F_{*}[0]$ if and only if
$$
\begin{cases}
W=0\\
f=\mbox{constant}\,.
\end{cases}
$$
It follows that $\ker F_*[0]=$Span$_{\R}(R_{\alpha})\subset
C^{\infty}(\mathcal{N}(L))$. The implicit function theorem between
Banach spaces implies that the Moduli space
$\mathfrak{M}(L)$ is a 1-dimensional smooth manifold.
\end{proof}
\begin{rem}\emph{
Note that the dimension of $\mathfrak{M}(L)$ does not depend on that
one of $L$. This is quite different from the Calabi-Yau case, where
the dimension of the Moduli space of deformations of special
Lagrangian submanifolds near a fixed compact $L$ is equal to the
first Betti number of $L$. This difference can be explained in the
following way: the deformations parameterized by curves tangent to
the contact structure are trivial, while those one along the Reeb
vector field $R_{\alpha}$ parameterize the Moduli space.}
\end{rem}
Now we study the following\\
\newline
\textbf{Extension problem:} Let $(M,\alpha_t,J_t,\epsilon_t)$,
$t\in(-\delta,\delta)$,
be a smooth family of contact Calabi-Yau manifolds. Given a compact
special Legendrian submanifold $p\colon L\hookrightarrow M$ of
$(M,\alpha_0,J_0,\epsilon_0)$ does it exist a family
$p_t\colon L\hookrightarrow M$ of special Legendrian submanifolds of $(M,\alpha_t,J_t,\epsilon_t)$
such that $p_0\colon L\hookrightarrow M$ coincides with $p$ ?\\
\newline
This is a  contact version of the extension problem in the Calabi-Yau
case (see \cite{LP} and  \cite{Pa}).\\
We can state the following
\begin{theorem}
\label{LP} Let $(M,\alpha_t,J_t,\epsilon_t)_{t\in (-\delta,\delta)}$
be a smooth family of contact Calabi-Yau manifolds. Let $p\colon L
\hookrightarrow M$ be a compact special Legendrian submanifold of
$(M,\alpha_0,J_0,\epsilon_0)$. Then there exists, for small $t$, a
family of compact special Legendrian submanifolds $p_t\colon L
\hookrightarrow (M,\alpha_t,J_t,\epsilon_t)$ such that $p_0=p$ if
and only if the condition
\begin{equation}
\label{extension}
[p^*(\Im\mathfrak{m}\,\epsilon_t)]=0
\end{equation}
holds for $t$ small enough.
\end{theorem}
\begin{proof}
The condition \eqref{extension} is necessary. Indeed if we can extend
$L$, then
$\Im\mathfrak{m}\,\epsilon_t$ is a closed form such that
$p_t^*(\Im\mathfrak{m}\,\epsilon_t)=0$. Since $p_t$ is homotopic to
$p_0$ we have
$$
[p_0^*(\Im\mathfrak{m}\,\epsilon_t)]=[p_t^*(\Im\mathfrak{m}\,\epsilon_t)]=0\,.
$$
In order to prove that condition \eqref{extension} is sufficient, we
can consider the map
$$
G\colon (-\sigma,\sigma)\times C^{1,\alpha}(J(p_{*}TL)) \to
C^{0,\alpha}(\Lambda^2(L))\oplus C^{0,\alpha}(\Lambda^{n}(L))
$$
defined as
$$
G(t,Z)=(\exp_{Z}^*(d\alpha_t),2\exp^*_Z(\Im \mathfrak{m}\,\epsilon_t))\,.
$$
By our assumption it follows that
$$
\mbox{Im}(G)\subset
dC^{1,\alpha}(\Lambda^1(L))\oplus
dC^{(1,\alpha)}(\Lambda^{n-1}(L))\,.
$$
Let $X\in C^{1,\alpha}(p_*(TL))$; a direct
computation and lemma \ref{linear} give
$$
\begin{aligned}
G_{*}[(0,0)](0,JX)=&(dp^*(\iota_{JX}d\alpha_0),2dp^*(\iota_{JX}\Im \mathfrak{m}\,\epsilon))\\
                 =&(dp^*(\iota_{JX}d\alpha_0),-d*p^*(\iota_{JX}d\alpha_0))\,,
\end{aligned}
$$
where $*$ is the Hodge operator of the metric
$p^*(g_{J})$  with respect to the
volume form $p^*(\Re\mathfrak{e}\,\epsilon)$. It follows that
$G_{*}[(0,0)](0,\cdot)$ is surjective and that
$$
\ker G_{*}[(0,0)]_{\{0\}\times C^{1,\alpha}(p_{*}(J(TL)))}\equiv \mathcal{H}^1(L)\,,
$$
where $\mathcal{H}^1(L)$ denotes the space of harmonic $1$-forms on $L$.\newline
Let
$$
A=\{X\in C^{1,\alpha}(p_*(TL))\,\,|\,\,p^*(\iota_{JX}d\alpha)\in
dC^{1,\alpha}(L)\oplus d^*C^{1,\alpha}(\Lambda^2(L))\}
$$
and
$$
\hat{G}=G_{|(-\delta,\delta)\times A}\,.
$$
Then by the Hodge decomposition of $\Lambda (L)$ it follows that
$$
G_{*}[(0,0)]_{\{0\}\times A}\colon A\to dC^{1,\alpha}(L)\oplus d^*C^{1,\alpha}(\Lambda^2(L))
$$
is an isomorphism. Again by the implicit function theorem and the
elliptic regularity there exists a local smooth solution of the
equation
$$
\hat{G}(t,\psi(t))=0\,.
$$
The extension of $p\colon L\hookrightarrow M$ is obtained by considering
$$
p_{t}:=\exp_{\psi(t)}\,.
$$
\end{proof}
\section{The 5-dimensional nilpotent case}
In this section we study invariant contact Calabi-Yau structures on
5-dimensional nilmanifolds. We will prove that a compact 5-dimensional
nilmanifold carrying an invariant Calabi-Yau structure is covered by a
Lie group whose Lie algebra is isomorphic to
$$
\mathfrak{g}=(0,0,0,0,12+34)\,,
$$
just described in section 2. Notation $\mathfrak{g}=(0,0,0,0,12+34)$
means that there exists a basis $\{\alpha_1,\dots,\alpha_{5}\}$ of the
dual space of the Lie algebra $\mathfrak{g}$ such that
$$
d\alpha_1=d\alpha_2=d\alpha_3=d\alpha_4=0\,,\quad d\alpha_5=\alpha_{1}\wedge\alpha_2+\alpha_{3}\wedge\alpha_4\,.
$$

First of all we note that 5-dimensional contact Calabi-Yau manifolds are in
particular Hypo. Recall that an \emph{Hypo structure} on a 5-dimensional
manifold is the datum of $(\alpha,\omega_1,\omega_2,\omega_3)$, where
$\alpha\in\Lambda^1(M)$ and $\omega_{i}\in\Lambda^{2}(M)$ and
\begin{itemize}
\item[1.] $\omega_{i}\wedge\omega_j=\delta_{ij}\,v$, for some
  $v\in\Lambda^{4}(M)$ satisfying $v\wedge\alpha\neq 0$;
\vspace{0.3cm}
\item[2.] $\iota_X\omega_1=\iota_Y\omega_2\iff
  \omega_{3}(X,Y)\geqslant 0$:
\vspace{0.3cm}
\item[3.] $d\omega_1=0$, $d(\omega_2\wedge\alpha)=0$, $d(\omega_3\wedge\alpha)=0$.
\end{itemize}
These structures have been introduced and studied by D. Conti and S. Salamon in
\cite{CS}.
Let $(M,\alpha,J,\epsilon)$ be a contact Calabi-Yau manifold of
dimension 5. Then
$$
\alpha\,,\quad \omega_1=\frac{1}{2}d\alpha\,,\quad
\omega_2=\Re\mathfrak{e}\,\epsilon\,,\quad \omega_3=\Im\mathfrak{m}\,\epsilon\,,
$$
define an Hypo structure on $M$.\\

The following lemma, whose proof is immediate, will be useful in
the sequel
\begin{lemma}
Let $M=G/\Gamma$ be a nilmanifold of dimension $5$. If $M$ admits an
invariant contact form, then the Lie algebra of $G$ is isomorphic to
one of the following 
\begin{itemize}
\item $(0,0,12,13,14+23)\,;$
\vspace{0.3cm}
\item $(0,0,0,12,13+24)\,;$
\vspace{0.3cm}
\item $(0,0,0,0,12+34)\,.$
\end{itemize}
\end{lemma}
Let $\mathfrak{g}$ be a non-trivial 5-dimensional nilpotent Lie algebra
and denote by $V=\mathfrak{g}^*$ the dual vector space of
$\mathfrak{g}$. There exists a filtration of $V$
$$
V^1\subset V^2\subset V^3\subset V^4\subset V^5=V\,,
$$
with $dV^i\subset \Lambda^2 V^{i-1}$ and dim$_{\R}V^i=i$. We may choose
the filtration $V$ in such a way that $V^{2}\subset \ker d\subset
V^{4}$.

Let $(M=G/\Gamma,\alpha,\omega_1,\omega_2,\omega_3)$ be a
nilmanifold endowed with an invariant Hypo structure
$(\alpha,\omega_1,\omega_2,\omega_3)$
\begin{enumerate}
\item[1. ] Assume that $\alpha\in  V^{4}$. Then we have the following
  (see \cite{CS})
\begin{theorem}
If $\alpha\in V^4$, then $\mathfrak{g}$ is either $(0,0,0,0,12)$,
$(0,0,0,12,13)$, or $(0,0,12,13,14)$.
\end{theorem}
In particular if $(M,\alpha,J,\epsilon)$ is contact Calabi-Yau, then
$\alpha\in V^4$.
\vspace{0.3cm}
\item[2.] Assume that $\alpha\notin V^4$. We have (see again \cite{CS})
\begin{lemma}\label{lemma}
If $\alpha\notin V^4$ and all $\omega_i$ are closed, then $\alpha$ is
orthogonal to $V^4$.
\end{lemma}
\begin{theorem}\label{th}
If $\alpha$ is orthogonal to $V^4$, then $\mathfrak{g}$ is one of
$$
(0,0,0,0,12)\,,\quad (0,0,0,0,12+34)\,.
$$
\end{theorem}
\end{enumerate}
Let $(M,\alpha, J,\epsilon)$ be a contact Calabi-Yau manifold of
dimension 5 endowed with an invariant contact Calabi-Yau structure;
then by $1.$ $\alpha$ does not belong to $V^4$. By lemma \ref{lemma}
$\alpha$ is orthogonal to $V^4$ and by theorem \ref{th} $
\mathfrak{g}=(0,0,0,0,12+34)\,.$ Hence we have proved the following
\begin{theorem}
Let $M=G/\Gamma$ be a nilmanifold of dimension $5$ admitting an
invariant contact Calabi-Yau structure. Then $\mathfrak{g}$ is
isomorphic to
$$
(0,0,0,0,12+34)\,.
$$
\end{theorem}
\section{Calabi-Yau manifolds of codimension $r$.}
In this section we extend the definition  of contact Calabi-Yau
manifold to codimension $r$ showing the analogous of theorem
\ref{LP}.

Let us consider the following
\begin{definition}
Let $M$ be a $2n+r$-dimensional manifold. An $r$-\emph{contact
structure} on $M$ is the datum
$\mathcal{D}=\{\alpha_1,\dots,\alpha_r\}$, where
$\alpha_i\in\Lambda^1(M)$, such that
\begin{itemize}
\item $d\alpha_1=d\alpha_{2}=\dots=d\alpha_{r}$;
\vskip0.2cm
\item $\alpha_{1}\wedge\dots\wedge\alpha_{r}\wedge (d\alpha_1)^n\neq 0$ .
\end{itemize}
\end{definition}
Note that if $\mathcal{D}=\{\alpha_1,\dots,\alpha_r\}$ is an
$r$-contact structure and $\xi:=\bigcap\,\ker\alpha_i$, then
$(\xi,d\alpha_1)$ is a symplectic vector bundle on $M$ and there
exists a unique set of vector fields $\{R_1,\dots,R_{r}\}$
satisfying
$$
\alpha_{i}(R_{j})=\delta_{ij}\,,\quad\iota_{R_i}d\alpha_i=0\mbox{
for any }i,j=1,\dots,r\,.
$$
Let us denote by $\mathfrak{C}_{\kappa}(\xi)$ the set of complex
structures on $\xi$ calibrated by the symplectic form
$\kappa=\frac{1}{2}d\alpha_1$ and by $\Lambda^{r}_0(M)$ the set of $r$-forms
$\gamma$ on $M$ satisfying
$$
\iota_{R_{i}}\gamma=0\mbox{ for any } i=1,\dots,r.
$$
If $J\in\mathfrak{C}_{\kappa}(\xi)$, then we extend it to $TM$
by defining
$$
J(R_{i})=0\,.
$$
Note that such a $J$ satisfies
$$
J^2=-I+\sum_{i=1}^r\,\alpha_i\otimes R_{i}\,.
$$
Consequently, for any  $J\in\mathfrak{C}_{\kappa}(\xi)$, we have $J(\Lambda^{r}_0(M))\subset \Lambda^{r}_0M$ and a natural
splitting of $\Lambda^{r}_0(M)\otimes\C$ in
$$
\Lambda^{r}_0(M)\otimes\C=\bigoplus_{p+q=r}\Lambda^{p,q}_J(\xi)\,.
$$
We can give the following
\begin{definition}
An $r$-contact Calabi-Yau manifold is the datum of
$(M,\mathcal{D},J,\epsilon)$, where
\begin{itemize}
\item $M$ is a $2n+r$-dimensional manifold;
\vskip0.2cm
\item $\mathcal{D}=\{\alpha_1,\dots,\alpha_r\}$ is an $r$-contact
structure;
\vskip0.2cm
\item $J\in\mathfrak{C}_{\kappa}(\xi)$
\vskip0.2cm
\item $\epsilon\in\Lambda^{n,0}_J(\xi)$ satisfies
$$
\begin{cases}
\epsilon\wedge\overline{\epsilon}=c_{n}\,\kappa^n\\
d\epsilon=0\,.
\end{cases}
$$
\end{itemize}
\end{definition}
\begin{ex}\label{ex2ccy}
\emph{ Let $M=H(3)/\Gamma\times S^1$ be the Kodaira-Thurston
manifold, where $H(3)$ is the $3$-dimensional Heisenberg  group and
$\Gamma$ is the lattice of $H(3)$ of matrices with integers entries.
Let
$$
\begin{aligned}
&\alpha_1=-2dz+2xdy\,,\\
&\alpha_2=-2dz+2xdy+2dt\,.
\end{aligned}
$$
One easily gets
$$
d\alpha_1=d\alpha_2=2dx\wedge dy
$$
and that $\mathcal{D}=\{\alpha_1,\alpha_2\}$ is a $2$-contact
structure on $M$. Note that $\xi=\ker\alpha_1\cap\ker\alpha_2$ is
spanned by $\{X_1=\partial_{x},X_2=\partial_{y}+x\partial_{z}\}$.
Moreover the Reeb fields of $\mathcal{D}$ are
$$
\begin{aligned}
&R_1=-\frac{1}{2}\partial_z-\frac{1}{2}\partial_t\,,\\
&R_2=\frac{1}{2}\partial_t\,.
\end{aligned}
$$
Therefore $\Lambda_0^1(M)$ is generated by $\{dx,dy\}$. Let
$J\in\mbox{End}(\xi)$ be the complex structure given by
$$
J(X_1)=X_2\,,\quad J(X_2)=-X_1
$$
and let $\epsilon\in\Lambda^{1,0}_J(\xi)$ be the form
$$
\epsilon=dx+idy\,.
$$
Then $(M,\mathcal{D},J,\epsilon)$ is a $2$-contact Calabi-Yau
structure.}
\end{ex}
As in the contact Calabi-Yau case if $(M,\mathcal{D},J,\epsilon)$ is
an $r$-contact Calabi-Yau manifold, then the $n$-form
$\Omega=\Re\mathfrak{e}\,\epsilon$ is a calibration on $M$. Moreover
an $n$-dimensional submanifold $p\colon L\hookrightarrow M$ admits an
orientation making it calibrated by $\Omega$ if and only if
$$
\begin{aligned}
&p^*(\alpha_{i})=0\mbox{ for any }\alpha_i\in\mathcal{D}\,,\\
&p^*(\Im\mathfrak{m}\,\epsilon)=0\,.
\end{aligned}
$$
A submanifold satisfying these equations will be called
\emph{special Legendrian}.
\begin{ex}
\emph{ Let $(M,\mathcal{D},J,\epsilon)$ be the $2$-contact
Calabi-Yau structure described in example \ref{ex2ccy}. Then
$$
L:=\left\{[A]\in H(3)/\Gamma\,|\, A=\left(
\begin{array}{ccc}
1  &x &0\\
0  &1 &0\\
0  &0 &1 \\
\end{array}
\right),\,x\in\R \right\}\times \{q\}\simeq S^1
$$
is a compact special Legendrian submanifold for any $q\in S^1$.}
\end{ex}
\noindent The proof of next theorem is very similar to that one of
theorem $\ref{LP}$ and it is omitted.
\begin{theorem}
\label{LP2} Let $(M,\mathcal{D}_t,J_t,\epsilon_t)_{t\in
(-\delta,\delta)}$ be a smooth family of  $r$-contact Calabi-Yau
manifolds. Let $p\colon L \hookrightarrow M$ be a compact special
Legendrian submanifold of $(M,\mathcal{D}_0,J_0,\epsilon_0)$. Then
there exists, for small $t$, a family of compact special Legendrian
submanifolds $p_t\colon L \hookrightarrow
(M,\mathcal{D}_t,J_t,\epsilon_t)$ extending $p\colon
L\hookrightarrow M$
 if and only if the
condition
\begin{equation*}
[p^*(\Im\mathfrak{m}\,\epsilon_t)]=0
\end{equation*}
holds for $t$ small enough.
\end{theorem}

\end{document}